\documentclass[a4paper,12pt]{amsart}
\usepackage{amssymb}
\usepackage{url}
\usepackage[T1]{fontenc}
\RequirePackage{calrsfs}
\DeclareSymbolFont{rsfscript}{OMS}{rsfs}{m}{b}
\DeclareSymbolFontAlphabet{\mathrsfs}{rsfscript}
\usepackage[utopia,expert]{mathdesign}
\usepackage{palatino}

\input xypic
\xyoption{all}
\xyoption{arc}

    \def\CM{{\mathbb{C}}}

    \def\FM{{\mathbb{F}}}

\def\IG{{\mathfrak I}}

    \def\QM{{\mathbb{Q}}}

    \def\ZM{{\mathbb{Z}}}

\def\Gb{{\mathbf G}}

\def\Lb{{\mathbf L}}

\def\Pb{{\mathbf P}}    \def\PC{{\mathcal{P}}}

\def\Sb{{\mathbf S}}    
    
\def\Ub{{\mathbf U}}

\def\Arm{{\mathrm{A}}}

\def\Drm{{\mathrm{D}}}

\def\Irm{{\mathrm{I}}}

\def\Rhat{{\hat{R}}}

\def\d{\delta}

\def\z{\zeta}

\def\mub{{\boldsymbol{\mu}}}

\DeclareMathOperator{\Hom}{{\mathrm{Hom}}}

\DeclareMathOperator{\Res}{{\mathrm{Res}}}

\def\to{\rightarrow}
\def\longto{\longrightarrow}

\def\fonctio#1#2#3#4{\begin{array}{ccc}
{#1} & \longto & {#2} \\
{#3} & \longmapsto & {#4} 
\end{array}}

\def\DS{\displaystyle}

\def\lexp#1#2{\kern\scriptspace\vphantom{#2}^{#1}\kern-\scriptspace#2}

\def\ge{\hspace{0.1em}\mathop{\geqslant}\nolimits\hspace{0.1em}}

\mathchardef\inferieur="321E
\mathchardef\superieur="321F

\def\eqna{\begin{eqnarray*}}
\def\endeqna{\end{eqnarray*}}

\def\gfq{{\FM_{\! q}}}

\def\gfqplus{{\FM_{\! q}^+}}

\catcode`\@=11
\long\def\@car#1#2\@nil{#1}
\long\def\@first#1#2{#1}
\long\def\@second#1#2{#2}
\long\def\ifempty#1{\expandafter\ifx\@car#1@\@nil @\@empty
  \expandafter\@first\else\expandafter\@second\fi}
\catcode`\@=12

\def\IGh{{\hat{\IG}}}

\begin{document}

\baselineskip=16pt

\title{A progenerator for representations of $\Sb\Lb_n(\gfq)$ \\ 
in transverse characteristic}

\author{C\'edric Bonnaf\'e}
\address{\noindent CNRS - UMR 5149 \\
Institut de Math\'ematiques et de Mod\'elisation de Montpellier \\
Universit\'e Montpellier 2 \\
Place Eug\`ene Bataillon \\
34095 MONTPELLIER Cedex \\
FRANCE} 

\makeatletter
\email{cedric.bonnafe@math.univ-montp2.fr}
\urladdr{http://ens.math.univ-montp2.fr/~bonnafe/}

\makeatother


\date{\today}

\begin{abstract} 
Let $G=\Gb\Lb_n(\gfq)$, $\Sb\Lb_n(\gfq)$ or $\Pb\Gb\Lb_n(\gfq)$, 
where $q$ is a power of some prime number $p$, 
let $U$ denote a Sylow $p$-subgroup of $G$ and let $R$ be a 
commutative ring in which 
$p$ is invertible. Let $\Drm(U)$ denote the derived subgroup of $U$ and 
let $e=\frac{1}{|\Drm(U)|} \sum_{u \in \Drm(U)} u$. The aim of this note is 
to prove that the $R$-algebras 
$RG$ and $eRGe$ are Morita equivalent (through the natural functor 
$RG$-mod $\longto$ $eRGe$-mod, $M \mapsto eM$).
\end{abstract}

\bigskip

\maketitle

\pagestyle{myheadings}

\markboth{\sc C. Bonnaf\'e}{\sc A progenerator for representations of $\Sb\Lb_n(\gfq)$}

Let $n$ be a non-zero natural number, 
$p$ a prime number, $q$ a power of $p$ and let $\gfq$ denote a finite field 
with $q$ elements. Let $G_n=\Sb\Lb_n(\gfq)$. We denote by $U_n$ 
the group of $n \times n$ unipotent upper 
triangular matrices with coefficients in $\gfq$ (so that $U_n$ is 
a Sylow $p$-subgroup of $G_n$). Let 
$\Drm(U_n)$ denote its derived subgroup: then, with $N=(n-1)(n-2)/2$, 
$$\Drm(U_n) = 
\Bigl\{
\begin{pmatrix}
1 & 0 & a_1 & \cdots & \cdots & a_{n-2} \\
0 & 1 & 0 & \ddots & & \vdots \\
\vdots & \ddots & \ddots & \ddots & \ddots & \vdots \\
\vdots & & \ddots & \ddots & \ddots & a_N\\
\vdots &  &  & \ddots & 1 & 0 \\
0 &  \cdots & \cdots & \cdots & 0 & 1 \\
\end{pmatrix}~\Bigl|~a_1,a_2,\dots,a_N \in \gfq\Bigr\}.$$
We fix a commutative ring $R$ in which $p$ is invertible and we set
$$e_n=\frac{1}{|\Drm(U_n)|} \sum_{u \in \Drm(U_n)} u\qquad \in \quad R\Drm(U_n).$$
Then $e_n$ is an idempotent of $RG_n$. The aim of this note 
is to prove the following result (recall that an idempotent $i$ of 
a ring $A$ is called {\it full} if $A=AiA$):

\bigskip

\def\module{{\mathrm{-mod}}}

\noindent{\bf Theorem 1.} 
{\it If $p$ is invertible in $R$, then $e_n$ is a full idempotent of $RG_n$.}

\bigskip

\begin{proof}
First, let $R_0=\ZM[1/p]$, let $\z$ be a primitive $p$-th root of unity in $\CM$ 
and let $\Rhat_0=R_0[\z]$. Let $\IG_0=R_0G_ne_nR_0G_n$ and $\IGh_0=\Rhat_0 G_ne_n \Rhat_0 G_n$. 
Since $p$ is invertible in $R$, there is a unique morphism of rings $R_0 \to R$ 
which extends to a morphism of rings $R_0 G_n \to RG_n$. So if $1 \in \IG_0$, 
then $1 \in \IG$. Also, as $(1,\z,\dots,\z^{p-2})$ is an $R_0$-basis of 
$\Rhat_0$, it is also an $R_0G_n$-basis of $\Rhat_0 G_n$. Therefore, 
if $1 \in \Rhat_0 G_n e_n \Rhat_0 G_n=\Rhat_0 \otimes_{R_0} (R_0G_n eR_0G_n)$, 
then $1 \in \IG_0$.
Consequently, in order to prove Theorem 1, we may (and we shall) work under the following 
hypothesis:

\medskip

\begin{quotation}
\noindent{\bf Hypothesis.} 
{\it From now on, and until the end of this proof, we assume that 
$R=\ZM[1/p,\z]$.}
\end{quotation}

\medskip

Now, let $P_n$ denote the subgroup of $\Sb\Lb_n(\gfq)$ defined by
$$P_n=\Bigl\{
\left(
\begin{array}{ccc|c}
&&& a_1 \\
& M && \vdots \\
&&& a_{n-1} \\
\hline
\vphantom{A^{\DS{A}}}0 & \cdots & 0 & 1
\end{array}\right)~\Bigl|~M \in \Sb\Lb_{n-1}(\gfq)
\quad\text{and}\quad a_1,\dots,a_{n-1} \in \gfq\Bigr\}.
$$
Then $U_n \subset P_n$. 
We shall prove by induction on $n$ that 
$$\text{\it $e_n$ is a full idempotent of $RP_n$}.\leqno{(\PC_n)}$$
It is clear that Theorem 1 follows immediately from $(\PC_n)$. 

\medskip

As $e_1=1$ and $e_2=1$, it follows that $(\PC_1)$ and $(\PC_2)$ hold. 
So assume that $n \ge 3$ and that $(\PC_{n-1})$ holds. 
Let $\Irm_n$ denote the identity $n \times n$ matrix and let 
$$V_n=\Bigl\{
\left(
\begin{array}{ccc|c}
&&& a_1 \\
& \Irm_{n-1} && \vdots \\
&&& a_{n-1} \\
\hline
\vphantom{A^{\DS{A}}}0 & \cdots & 0 & 1
\end{array}\right)~\Bigl|~a_1,\dots,a_{n-1} \in \gfq\Bigr\}.
$$
Then $V_n \simeq (\gfqplus)^{n-1}$ and 
$P_n = \Sb\Lb_{n-1}(\gfq) \ltimes V_n \simeq \Sb\Lb_{n-1}(\gfq) \ltimes (\gfqplus)^{n-1}$. 
We set $V_n'=\Drm(U_n) \cap V_n$, so that $V_n' \simeq (\gfqplus)^{n-2}$ is normalized 
by $P_{n-1}$. Then 
$$\Drm(U_n) = \Drm(U_{n-1}) \ltimes V_n'.$$
We now define
$$f_n=\frac{1}{|V_n'|} \sum_{v \in V_n'} v,$$
so that 
$$e_n = e_{n-1} f_n.$$
By the induction hypothesis, there exists $g_1$, $h_1$,\dots, $g_l$, $h_l$ in 
$P_{n-1}$ and $r_1$,\dots, $r_l$ in $R$ such that 
$$1=\sum_{i=1}^l r_i g_i e_{n-1} h_i.$$
Therefore, as $P_{n-1}$ normalizes $V_n'$, it centralizes $f_n$ and so 
$$f_n = \Bigl(\sum_{i=1}^l r_i g_i e_{n-1}  h_i \Bigr) f_n
= \sum_{i=1}^l r_i g_i e_{n-1} f_n h_i
= \sum_{i=1}^l r_i g_i e_n h_i.$$
So $f_n \in RP_n e_n RP_n$. 

Let $\mub_p$ denote the subgroup of $R^\times$ generated by $\z$. 
If $\chi \in \Hom(V_n,\mub_p)$, we define $b_\chi$ to be the associated primitive 
idempotent of $RV_n$:
$$b_\chi = \frac{1}{|V_n|} \sum_{v \in V_n} \chi(v)^{-1} v \qquad \in \quad RV_n.$$
Then, as $V_n$ is an elementary abelian $p$-group, we get
$$f_n = \sum_{\substack{\chi \in \Hom(V_n,\mub_p) \\ \Res_{V_n'}^{V_n} \chi = 1}} b_\chi.$$
We fix a non-trivial element $\chi_0 \in \Hom(V_n,\mub_p)$ whose restriction 
to $V_n'$ is trivial. Then 
$$b_{\chi_0} = b_{\chi_0} f_n\qquad\text{and}\qquad b_1 = b_1 f_n,$$
so $b_1$ and $b_{\chi_0}$ belong to $RP_n e_n RP_n$. 

But $\Sb\Lb_{n-1}(\gfq) \subset P_n$ has only two orbits for its action on 
$\Hom(V_n,\mub_p)$: the orbit of $1$ and the orbit of $\chi_0$. Therefore, 
$b_\chi \in RP_n e_n RP_n$ for all $\chi \in \Hom(V_n,\mub_p)$. 
Consequently,
$$1=\sum_{\chi \in \Hom(V_n,\mub_p)} \!b_\chi ~\in~ RP_n e_n RP_n,$$
as desired.
\end{proof}

\bigskip

\noindent{\bf Finite reductive groups.} 
Let $\FM$ be an algebraic closure of $\gfq$, let $\Gb$ be a connected 
reductive group over $\FM$ and let $F : \Gb \to \Gb$ be an isogeny 
such that some power $F^\d$ is a Frobenius endomorphism relative to 
an $\gfq$-structure. We denote by $\Ub$ an $F$-stable maximal unipotent 
subgroup of $\Gb$ (it is the unipotent radical of an $F$-stable Borel subgroup). 
Define 
$$e=\frac{1}{|\Drm(\Ub)^F|} \sum_{u \in \Drm(\Ub)^F} u\qquad \in \quad R\Gb^F.$$
The next result follows immediately from Theorem 1:

\bigskip

\noindent{\bf Theorem 2.} 
{\it Assume that $(\Gb,F)$ is split of type $\Arm$. Then $e$ is a full 
idempotent of $R\Gb^F$.}

\bigskip

%

\noindent{\bf Corollary 3.} 
{\it If $(\Gb,F)$ is split of type $\Arm$, then the functors 
$$\fonctio{R\Gb^F\module}{e R\Gb^F e\module}{M}{eM}
\qquad\text{\it and}\qquad
\fonctio{R\Gb^F e\module}{R\Gb^F\module}{N}{R\Gb^F e \otimes_{eR\Gb^F e} N}$$
are mutually inverse equivalences of categories. In particular, 
$R\Gb^F$ and $eR\Gb^F e$ are Morita equivalent, and $R\Gb^F e$ is a progenerator 
for $R\Gb^F$.}

\begin{proof}
This follows from Theorem 2 and, for instance, \cite[Example 18.30]{lam}.
\end{proof}

\bigskip

\noindent{\bf Examples.} 
Theorem 2 and Corollary 3 can be applied for instance in the case where 
$\Gb^F = \Gb\Lb_n(\gfq)$, $\Sb\Lb_n(\gfq)$ or $\Pb\Gb\Lb_n(\gfq)$.

\bigskip

\noindent{\bf Comments.} 
(1) It is natural to ask whether Theorem 2 (or Corollary 3) can be generalized 
to other finite reductive groups. In fact, it cannot be generalized: 
indeed, if for instance $R=\CM$, then saying that $e$ is a full idempotent 
of $R\Gb^F$ means that every irreducible character of $\Gb^F$ is an irreducible 
component of an Harish-Chandra induced of some Gelfand-Graev character. 
But, if $\Gb$ is quasi-simple and $(\Gb,F)$ is not split of type $A$, 
then $\Gb^F$ contains a unipotent 
character which does not belong to the principal series: this character 
cannot be an irreducible component of an Harish-Chandra induced of 
a Gelfand-Graev character.

\medskip

(2) In \cite{BR}, a crucial step for the proof of a special case of the geometric 
version of Brou\'e's abelian defect conjecture was \cite[Theorem 4.1]{BR}, 
where R. Rouquier and the author have proved the above Theorem 2 in the case 
where $R$ is the integral closure of $\ZM_{\! \ell}$) in a sufficiently large 
algebraic extension of $\QM_{\! \ell}$ (here, $\ell$ is a prime number different from $p$). 
The proof was essentially based on the classification, 
due to Dipper \cite[4.15 and 5.23]{dipper}, of simple modules for $G_n$ 
in characteristic $\ell$, and especially of cuspidal ones, which involves 
Deligne-Lusztig theory.

The interest of the proof given here is that it does not rely 
on any classification of simple modules, and is based on elementary methods: 
as a by-product of this elementariness, Theorem 2 and Corollary 3 are valid 
over any commutative ring (in which $p$ is invertible, which is a necessary condition if one 
wants the idempotent $e_n$ to be well-defined).

\bigskip

%


\begin{thebibliography}{ABC}
\bibitem[1]{BR} {\sc C. Bonnaf\'e \& R. Rouquier}, 
Coxeter orbits and modular representations, 
{\it Nagoya Math. J.} {\bf 183} (2006), 1-34.

\medskip

\bibitem[2]{dipper} {\sc R. Dipper}, 
On quotients of Hom-functors and representations of finite general
linear groups II, {\it J. Algebra} {\bf 209} (1998), 199-269.

\medskip

\bibitem[3]{lam} {\sc T.-Y. Lam}, {\it Lectures on Modules and Rings}, 
Graduate Texts in Mathematics {\bf 189}, Springer, 1999, xxiv + 557 pages.
\end{thebibliography}
\end{document}